\documentclass[12pt,a4paper,reqno]{amsart}

\usepackage{anysize}
\usepackage[utf8]{inputenc}
\usepackage{amsmath}
\usepackage{amssymb}
\usepackage{array}

\begin{document}

\title{Dougall's $_5F_4$ sum and the WZ-algorithm}

\author{Jesús Guillera}
\address{Department of Mathematics, University of Zaragoza, 50009 Zaragoza, SPAIN}
\email{jguillera@gmail.com}

\date{}

\maketitle

\begin{abstract}
We show how to prove the examples of a paper by Chu and Zhang using the WZ-algorithm. 
\end{abstract}

\subsection*{Keywords}
Generalized hypergeometric series; Dougall's sum; WZ-algorithm; Ramanujan-type series; Evaluation of constants.
\subsection*{Mathematics Subject Classification}
33C20; 33F10; 11Y60

\section{Introduction}
Chu and Zhang had the good idea of accelerating the following Dougall's sum \cite{Chu-Zhang}, which is convergent for ${\rm Re}(1+2a-b-c-d-e)>0$:
\[
\sum_{k=0}^{\infty} \frac{(b)_k(c)_k(d)_k(e)_k}{(1+a-b)_k(1+a-c)_k(1+a-d)_k(1+a-e)_k}(a+2k),
\]
using many recurrence relations patterns. Their paper has interesting formulas and many nice examples. We show that the Wilf-Zeilberger algorithm can do it automatically \cite{PWZ}, and give four examples corresponding to Examples $42$, $45$, $50$ and $53$ of that paper. Note that Dougall's sum  in \cite{Chu-Zhang} has two wrong signs and observe that it is symmetric in $b$, $c$, $d$ and $e$ (that is, these parameters are interchangeable). 
\\ \par We will use the following notation:
\[
A(a,b,c,d,e,k)=\frac{(b)_k(c)_k(d)_k(e)_k}{(1+a-b)_k(1+a-c)_k(1+a-d)_k(1+a-e)_k}(a+2k).
\]
Then Dougall's sum is 
\[
\Omega(a,b,c,d,e)=\sum_{k=0}^{\infty} A(a,b,c,d,e,k).
\]
If $e \to -\infty$ then we have the sum
\[
\Phi(a,b,c,d)=\sum_{k=0}^{\infty} B(a,b,c,d,k) 
\]
where
\[
B(a,b,c,d,k)=(-1)^k \frac{(b)_k(c)_k(d)_k}{(1+a-b)_k(1+a-c)_k(1+a-d)_k}(a+2k).
\]
This introduction is short but the interested reader can also see the introduction and the references of \cite{Chu-Zhang}. 

\section{The WZ-algorithm as a black box}

Let $F(n,k)$ be hypergeometric in $n$ and $k$, that is a function such that $F(n+1,k)/F(n,k)$ and $F(n,k+1)/F(n,k)$ are rational functions. Then, we can use the Zeilberger's Maple package $\texttt{SumTools[Hypergeometric]);}$.
The output of $\texttt{Zeilberger(F(n,k),n,k,N)[1];}$ is an operator $O(N)$ of the following form
\[
O(N) = P_0(n)+P_1(n) \, N+P_2(n) \, N^2 + \cdots + P_{m}(n) \, N^{m},
\]
where $P_0(n), \, P_1(n), \dots, P_m(n)$ are polynomials of $n$, and $N$ is an operator which increases $n$ in $1$ unity, that is $N F(n,k)= F(n+1,k)$. The output of $\texttt{Zeilberger(F(n,k),n,k,N)[2];}$ gives a function $G(n,k)$ such that
\[
O(N)  F(n,k) = G(n,k+1)-G(n,k).
\]
If the operator $O(N)$ is just $N-1$, then one has
\begin{equation}\label{wz-pair}
F(n+1,k) - F(n,k) = G(n,k+1) - G(n,k),
\end{equation}
and we say that $(F,G)$ is a WZ-pair (Wilf-Zeilberger). H. Wilf and D. Zeilberger proved in a seminal paper \cite{Wilf-Zeil} (for which they received the Steele Prize), that in this case there exists a rational function $R(n,k)$ (the so called certificate), such that $G(n,k)=R(n,k) F(n,k)$; relation that is valid for all values of $n$ and $k$ once we have simplified the product. When the operator is not exactly $N-1$ but is of degree $1$, then is very easy to transform the function into another one with output equal to $N-1$. The interesting fact is that this is always the case for the functions $U(n,k)$ constructed from the functions inside the Dougall's sum in the way we show in the examples. That is, with independence of the recurrence pattern that one is taking, the operator has always degree equal to $1$. 
\par WZ-pairs are gems because of their interesting properties. For example, H. Wilf and D. Zeilberger, summing both sides of (\ref{wz-pair}) for $n \geq 0$, proved that
\begin{equation}\label{pro-1-wz}
\sum_{n=0}^{\infty} G(n,k+1) = \sum_{n=0}^{\infty} G(n,k) - F(0,k) + \lim_{n \to \infty} F(n,k),
\end{equation}
in case that the summation is convergent. If the above limit is equal to $0$, one has
\[
\sum_{n=0}^{\infty} G(n,k) = F(0,k)+\sum_{n=0}^{\infty} G(n,k+1),
\]
and from it, we get the following list of identities:
\begin{align*}
\sum_{n=0}^{\infty} G(n,0) &= F(0,0)+\sum_{n=0}^{\infty} G(n,1), \\
\sum_{n=0}^{\infty} G(n,1) &= F(0,1)+\sum_{n=0}^{\infty} G(n,2), \\
\sum_{n=0}^{\infty} G(n,2) &= F(0,2)+\sum_{n=0}^{\infty} G(n,3), \\
& \,  \, \, \vdots
\end{align*}
Then, from this infinite list, we get
\begin{equation}\label{pro-2-wz}
\sum_{n=0}^{\infty} G(n,0)=\sum_{k=0}^{\infty} F(0,k) + \lim_{k \to \infty}  \sum_{n=0}^{\infty} G(n,k),
\end{equation}
when both sides are convergent. If the convergence of the series is uniform then we can go further and commute the limit with the series.

\section{Examples} 
In this section we prove the Examples $42$, $45$, $50$ and $53$ of the paper \cite{Chu-Zhang} using methods based on the WZ-algorithm, but we can derive all the other examples of that paper using the same kind of WZ-techniques.

\subsection*{Example 42} 
It  accelerates 
\[
\Omega \left(\frac32,\frac12,1,1,1 \right)=\frac74 \zeta(3),
\]
see \cite[Theorem 5]{Chu-Zhang}, using the recurrence relation [3,0,1,2,2], see \cite[Lemma 22]{Chu-Zhang}. To do it with the WZ-method we consider the function
\[
U(n,k)=A\left( \frac32+3n, \frac12+0n, 1+1n,1+2n,1+2n,k \right),
\]
and use the Maple package \texttt{with(SumTools[Hypergeometric])}. We see that the output of  \texttt{Zeilberger(U(n,k),n,k,N)[1]} is equal to
\[
-1024 \left(\frac12+n\right)\left(1+n\right)^5 N+1728 \left(\frac32+n\right)^2 \left(\frac23+n\right) \left(\frac43+n\right)\left(\frac34+n\right)\left(\frac54+n\right).
\]
This output indicates that if
\begin{align*}
F(n,k) &=U(n,k) {\displaystyle \frac{\left(\frac12\right)_n(1)_n^5}{\left(\frac32\right)_n^2\left(\frac23\right)_n\left(\frac43\right)_n\left(\frac34\right)_n\left(\frac54\right)_n}} \left( - \frac{16}{27} \right)^n \\ &= \frac{U(n,k)}{2n+1} \, {\displaystyle \frac{(1)_n^5}{\left(\frac32\right)_n\left(\frac23\right)_n\left(\frac43\right)_n\left(\frac34\right)_n\left(\frac54\right)_n}} \left( - \frac{16}{27} \right)^n,
\end{align*}
then the output of \texttt{Zeilberger(F(n,k),n,k,N)[1]} is equal to $N-1$ (Observe that $-1024/1728=-16/27$). Hence, there is a function $G(n,k)$, such that $(F,G)$ is a WZ-pair. Writing
\begin{verbatim}
               simplify(Zeilberger(F(n,k),n,k,N)[2]/F(n,k))
\end{verbatim}
we get a rational function $R(n,k)$ such that $G(n,k)=R(n,k) F(n,k)$. In addition, we see that
\[
R(n,0)=\frac{172n^2+269n+106}{9(2+3n)(3+4n)}, \quad U(n,0)=\frac32+3n.
\]
As the limit of $F(n,k)$ as $n \to \infty$ is equal to $0$, we can use the identity (\ref{pro-2-wz}):
\[
\sum_{n=0}^{\infty} G(n,0) = \sum_{k=0}^{\infty} F(0,k) + \lim_{k \to \infty} \sum_{n=0}^{\infty} G(n,k).
\]
We can easily check that in this example the limit commutes with the summation and is equal to zero. Therefore, we have
\[
\sum_{n=0}^{\infty} G(n,0) = \sum_{k=0}^{\infty} F(0,k)=\Omega \left(\frac32,\frac12,1,1,1 \right)=\frac74 \zeta(3).
\]
Hence
\[
36 \sum_{n=0}^{\infty} G(n,0) = \sum_{n=0}^{\infty} {\displaystyle \frac{(1)_n^5}{\left(\frac32\right)_n\left(\frac53\right)_n\left(\frac43\right)_n\left(\frac74\right)_n\left(\frac54\right)_n}} \left( - \frac{16}{27} \right)^n (172n^2+269n+106)=63\zeta(3),
\]
which is the formula of \cite[Example 42]{Chu-Zhang}.

\subsection*{Example 45}
This example uses the same recurrence relation, that is [3,0,1,2,2], see \cite[Lemma 22]{Chu-Zhang} and our final remark. We consider the function
\[
U(n,k)=A\left( \frac12+3n, \frac12+0n, \frac12+1n,\frac12+2n,\frac12+2n,k \right).
\]
Then, proceeding as in the previous example, we obtain
\[
F(n,k)=U(n,k)\cdot n {\displaystyle \frac{\left(\frac12\right)_n\left(\frac14\right)_n^3\left(\frac34\right)_n^3}{\left(\frac13\right)_n\left(\frac23\right)_n\left(1\right)_n^5}} \left( - \frac{16}{27} \right)^n,
\]
and the rational function $R(n,k)$ such that $G(n,k)=R(n,k) F(n,k)$. In addition, we get
\[
R(n,0)=\frac{1376n^4+1808n^3+784n^2+138n+9}{16n(1+6n)(2+3n)(1+3n)}, \quad U(n,0)=\frac12+3n.
\]
We see that
\[
\sum_{k=0}^{\infty} F(0,k) = \sum_{k=0}^{\infty} \left\{ 0 \cdot A \left( \frac12, \frac12, \frac12, \frac12, \frac12, k  \right) \right\}= 0, \quad \text{and} \quad \Omega \left( \frac12, \frac12, \frac12, \frac12, \frac12 \right) = \infty.
\]
Therefore, see (\ref{pro-2-wz}):
\[
\sum_{n=0}^{\infty} G(n,0) = \sum_{k=0}^{\infty} F(0,k) + \lim_{k \to \infty} \sum_{n=0}^{\infty} G(n,k) =  \lim_{k \to \infty} \sum_{n=0}^{\infty} G(n,k),
\]
and we can check that the limit commutes with the summation, and in addition the limit of $G(n,k)$ as $k \to \infty$ is zero for $n \neq 0$. Hence
\[
\sum_{n=0}^{\infty} G(n,0) = \lim_{k \to \infty} G(0,k) = \lim_{k \to \infty} \left( \frac{1}{32} \frac{\left(\frac12\right)_k^4}{(1)_k^4} \frac{(2k+3)^2(2k+1)^3}{(k+1)^2(k+2)} \right) = \frac{1}{\pi^2}.
\]
That is
\[
\sum_{n=0}^{\infty} {\displaystyle \frac{\left(\frac12\right)_n\left(\frac14\right)_n^3\left(\frac34\right)_n^3}{\left(\frac43\right)_n\left(\frac53\right)_n\left(1\right)_n^5}} \left( - \frac{16}{27} \right)^n (1376n^4+1808n^3+784n^2+138n+9)=\frac{64}{\pi^2},
\]
which is the formula of \cite[Example 45]{Chu-Zhang}. We will finally deduce that this example accelerates \[
9-\frac{9}{16} \, \Omega \left( \frac72, \frac12, \frac32, \frac52, \frac52 \right).
\]
For it observe that if $(F(n,k),G(n,k))$ is a WZ-pair then $(F(n+1,k),G(n+1,k)$ is also a WZ-pair. Hence
\[
\sum_{n=0}^{\infty}  G(n+1,0) = \sum_{k=0}^{\infty} F(1,k) + \lim_{k \to \infty} \sum_{n=0}^{\infty} G(n+1,k).
\]
As the above limit is equal to zero, we obtain
\begin{align*}
\frac{1}{\pi^2} = \sum_{n=0}^{\infty}  G(n,0) &= G(0,0)+\sum_{n=0}^{\infty}  G(n+1,0) = G(0,0)+\sum_{k=0}^{\infty} F(1,k) \\ &= \frac{9}{64}-\frac{9}{1024} \, \Omega \left( \frac72, \frac12, \frac32, \frac52, \frac52 \right),
\end{align*}
and from it, we deduce that
\[
\Omega \left( \frac72, \frac12, \frac32, \frac52, \frac52 \right) = 16 - \frac{1024}{9\pi^2},
\]
which is a known result because Maple also gives the exact value of this Dougall's sum.

\subsection*{Example 50} 
It  accelerates 
\[
\Phi \left(\frac32,1,1,1 \right)=G \quad \text{(Catalan constant)},
\]
see \cite[Theorem 5]{Chu-Zhang}, using the recurrence relation [3,0,2,2], see \cite[Lemma 22]{Chu-Zhang}. To do it with the WZ-method we consider the function
\[
U(n,k)=B \left( \frac32+3n, \frac12+0n,1+2n,1+2n,k \right),
\]
Then, with the Wilf-Zeilberger's algorithm, we obtain
\[
F(n,k) = U(n,k) {\displaystyle \frac{\left(\frac12\right)_n^2(1)_n^2}{\left(\frac32\right)_n^2\left(\frac56\right)_n\left(\frac76\right)_n}} \left( \frac{16}{27} \right)^n = \frac{U(n,k)}{(2n+1)^2} \, {\displaystyle \frac{(1)_n^2}{\left(\frac56\right)_n\left(\frac76\right)_n}} \left(\frac{16}{27} \right)^n,
\]
and the rational function $R(n,k)$ such that $G(n,k)=R(n,k) F(n,k)$. In addition, we get
\[
R(n,0)=\frac{22n+21}{5+6n}, \quad U(n,0)=\frac12+3n,
\]
and, see (\ref{pro-2-wz}):
\[
\sum_{n=0}^{\infty} G(n,0) = \sum_{k=0}^{\infty} F(0,k) = \Phi \left(\frac32,1,2,2 \right) = G.
\]
Hence
\[
30 \sum_{n=0}^{\infty} G(n,0)= \sum_{n=0}^{\infty} {\displaystyle \frac{(1)_n^2}{\left(\frac76\right)_n\left(\frac{11}{6}\right)_n}} \frac{22n+21}{2n+1} \left( \frac{16}{27} \right)^n = 30 \, G,
\]
which is the formula in \cite[Example 50]{Chu-Zhang}.

\subsection*{Example 53}
This example uses the recurrence relation [3,1,1,1,3], see \cite[Lemma 25]{Chu-Zhang}. From the function
\[
U(n,k)=A\left( \frac12+3n, \frac12+n, \frac12+n,\frac12+n,3n,k \right),
\]
we obtain
\[
F(n,k)=U(n,k)\cdot n^2 {\displaystyle \frac{\left(\frac12\right)_n^3 \left(\frac13\right)_n\left(\frac23\right)_n}{\left(1\right)_n^5}} \left(\frac{27}{64} \right)^n,
\]
and the rational function $R(n,k)$ such that $G(n,k)=R(n,k) F(n,k)$. In addition, we get
\[
R(n,0)=\frac{74n^2+27n+3}{16n^2(1+6n)}, \quad U(n,0)=\frac12+3n,
\]
and, see (\ref{pro-2-wz}):
\[
G(0,k)=\frac{3}{32} \frac{\left(\frac12\right)_k^2 }{(1)_k^2} \, \frac{(2k+1)^2}{(k+1)^2}.
\]
Then, we use the Wilf-Zeilberger's identity (\ref{pro-1-wz}):
\[
\sum_{n=0}^{\infty}  G(n,k+1) = \sum_{n=0}^{\infty} G(n,k) - F(0,k) + \lim_{n \to \infty} F(n,k).
\]
But $F(0,k)$ and the limit of $F(n,k)$ as $n \to \infty$ are equal to zero in this example. Therefore, we have
\[
\sum_{n=0}^{\infty}  G(n,k+1) = \sum_{n=0}^{\infty} G(n,k).
\]
In addition, we can prove that Carlson's Theorem holds here. Hence
\[
\sum_{n=0}^{\infty}  G(n,k) = {\rm constant} \quad \forall k \in \mathbb{C}.
\]
To determine the constant, we write
\[
\sum_{n=0}^{\infty}  G(n,k) = \lim_{k \to -\frac12} \sum_{n=0}^{\infty} G(n,k) = \lim_{k \to -\frac12} G(0,k) = \frac{3}{2\pi^2},
\]
because the limit commutes with the summation and is equal to $0$ when $n \neq 0$. Hence
\[
32 \sum_{n=0}^{\infty} G(n,0)=\sum_{n=0}^{\infty} {\displaystyle \frac{\left(\frac12\right)_n^3 \left(\frac13\right)_n \left(\frac23\right)_n}{\left(1\right)_n^5}} \left( \frac{27}{64} \right)^n (74n^2+27n+3)=\frac{48}{\pi^2},
\]
which is the formula \cite[Example 53]{Chu-Zhang}.

\subsection*{Remark}\label{obser}
Some of the examples in \cite{Chu-Zhang}, like the Example $45$, have $1+2a-b-c-d-e \leq 0$, and then $\Omega(a,b,c,d,e)$ is not defined. Hence the left hand side of the formulas in the theorems of \cite{Chu-Zhang} has no sense in these cases and variants of the theorems would be needed. In this paper we have shown that the Wilf-Zeilberger's techniques can manage all the situations.

\section*{ADDENDUM: A MAPLE PROGRAM (May 20, 2019)} \noindent
The following sections are new and do not appear in the reference journal.

\section{How to use the program}

Dougall's sums are $\Omega(a,b,c,d,e)$ and $\Omega(a,b,c,d,\infty)=\Phi(a,b,c,d)$.
In next section we provide a Maple program (written in Maple 9) which finds WZ pairs corresponding to the application of the kind of recurrences explained before. It does it in an automatic way with the help of the very powerful procedure \texttt{Zeilberger}. Using the program is very easy. First write an \texttt{entry}; then \texttt{SER(entry);} gives a {\it guessing} of the evaluation in close form of the corresponding series. Then, \texttt{FF(entry);} gives a function $F(n,k)$ and \texttt{CC(entry);} a rational function $R(n,k)$ that certificates that the pair formed by $F(n,k)$ and $G(n,k)=R(n,k) F(n,k)$ is a Wilf-Zeilberger pair. When $e=\infty$ use the notation oo for the last value of the entry. Here are three examples:

\begin{verbatim}

entry:=1+4*n,1/2+2*n,1/2+1*n,1/2+1*n,oo; SER(entry); FF(entry); CC(entry);

entry:=3/2+2*n,1+0*n,1/2+1*n,1+1/2*n,1/2+1/2*n; 

entry:=1/2+2*n,1/4+1/2*n,3/4+1/2*n,1/2+1*n,1/2+0*n; 

\end{verbatim}

The numbers in the recurrences used in \cite{Chu-Zhang-bis} and \cite{Chu} are positive integers but there are also some other possibilities (see the two last examples above).
The reader can try other entries; e.g. taken from the examples given by Chu and Zhang in \cite{Chu-Zhang-bis}, and by Chu in \cite{Chu}.

\section{ALMOST AUTOMATIC PROOFS BY COMPUTER (a Maple program)}

\begin{verbatim}

restart; with(SumTools[Hypergeometric]):

po:=proc(i,j) local f;
f:=i->i-floor(i): if type(i,numeric) then 
if i<=0 and f(i)=0 then return pochhammer(1,j)/product(j-h,h=0..-i): 
else if i>0 and f(i)=0 then return
simplify(pochhammer(i,j)/pochhammer(f(i)+1,j))*pochhammer(f(i)+1,j) 
else: if f(i) > 0 then return
simplify(pochhammer(i,j)/pochhammer(f(i),j))*pochhammer(f(i),j) 
fi: fi: fi: else return pochhammer(i,j): fi: end:

UU:=proc(r1,r2,r3,r4,r5) global U,po; local en,qq,ne,A1,A2,A; 
alias(p=pochhammer): en:=[r1,r2,r3,r4,r5]:
if r5=oo then ne:=1 else ne:=0 fi:
A1:=(a,b,c,d,e,k)->(-1)^(ne*k)*product(po(en[i],k),i=2..nops(en)-ne):
A2:=(a,b,c,d,e,k)->product(po(1+en[1]-en[i],k),i=2..nops(en)-ne):
A:=(a,b,c,d,e,k)->A1(a,b,c,d,e,k)/A2(a,b,c,d,e,k)*(a+2*k):
qq:=r1,r2,r3,r4,r5: U:=(n,k)->A(qq,k): end:

FF:=proc(r1,r2,r3,r4,r5) global z0,F; local o,p1,p2,t1,t2,T,pz,z1; 
UU(r1,r2,r3,r4,r5): o:=Zeilberger(U(n,k),n,k,N)[1]: 
p1:=-solve(factor(coeff(o,N,1)),n): p2:=-solve(factor(coeff(o,N,0)),n):
pz:=subs(n=1/n,-factor(coeff(o,N,1))/factor(coeff(o,N,0))):
z0:=subs(n=0,(simplify(pz))):
t1:=map(po,[p1],n): t2:=map(po,[p2],n):
T:=(n,k)->product(t1[i],i=1..nops(t1))/product(t2[i],i=1..nops(t2)): 
F:=(n,k)->U(n,k)*T(n,k)*z0^n: z1:=Zeilberger(F(n,k),n,k,N)[1]: 
if z1=N-1 or z1=1-N then return("F(n,k)"=F(n,k)) fi end:

CC:=proc(r1,r2,r3,r4,r5) local sF; global R,G,ss; 
FF(r1,r2,r3,r4,r5):
sF:=simplify(Zeilberger(F(n,k),n,k,N)[2]/F(n,k)):
R:=(nn,kk)->simplify(subs(n=nn,k=kk,sF)):
G:=(n,k)->F(n,k)*R(n,k): return("R(n,k)"=R(n,k)): end:

SER:=proc(r1,r2,r3,r4,r5) local w,ss,res; 
CC(r1,r2,r3,r4,r5): w:=factor(G(n,0)): Digits:=30:
ss:=evalf(Sum(w,n=0..infinity)): 
if type(ss,numeric) then res:=identify(ss,extension=[Catalan]): 
else res:="there is no sum from n=0": fi: 
if abs(z0)>1 then res:="it is not convergent" fi:
return(Sum(w,n=0..infinity)=res): end:
 
\end{verbatim}

\newpage

\section{Other examples}

For proving the following ``divergent" Ramanujan-like series for $1/\pi^2$ (dual of the Zeilberger--Amdeberham formula for $\zeta(3)$, see \cite{amde-zeil}):
\[
\sum_{n=0}^{\infty} \frac{\left(\frac12\right)_n^5}{(1)_n^5} (205n^2+160n+32) (-1024)^n \, \, `` \! = \! " \, \, \frac{16}{\pi^2},
\]
we can use \, \texttt{entry:=0+3*n,0+2*n,0+2*n,0+2*n,0+2*n;}
\\ \\
In \cite[Page 497, Example 42]{Chu-Zhang-bis}, Chu and Zhang proved the following formula
\[
\sum_{n=0}^{\infty} \frac{(1)_n^5}{\left(\frac32\right)_n\left(\frac43\right)_n\left(\frac53\right)_n\left(\frac54\right)_n\left(\frac74\right)_n} (172n^2+269n+106) \left(-\frac{16}{27}\right)^n = 69 \zeta(3),
\]
with the recurrences $(3,0,1,2,2)$ and the values $a=3/2$, $b=1/2$, $c=d=e=1$. For proving directly the corresponding ``divergent" Ramanujan-like series for $1/\pi^2$, namely
\[
\sum_{n=0}^{\infty} \frac{\left(\frac12\right)_n\left(\frac13\right)_n\left(\frac23\right)_n\left(\frac14\right)_n\left(\frac34\right)_n}{(1)_n^5} (172n^2+75n+9) \left(-\frac{27}{16}\right)^n \, \, `` \! = \! " \, \, \frac{48}{\pi^2}.
\]
one can use \, \texttt{entry:=1/2+3*n,1/2+2*n,1/2+n,1/2+n,0+3*n;}   
\\ \\  
In \cite[Page 525, Example 14]{Chu}, Chu proved the following formula:
\[
\sum_{n=0}^{\infty} \frac{(1)_n^5}{\left(\frac32\right)_n^3\left(\frac43\right)_n\left(\frac53\right)_n} \frac{28n^2+38n+13}{(-27)^n} = \frac{21}{2} \zeta(3),
\]
with the recurrences $(3,0,1,1,1)$ and the values $a=1, b=c=d=e=1/2$, and bisection. We can avoid bisection with \, \texttt{entry:=3/2+2*n,1/2+1/2*n,1+1/2*n,1+n,1+n;} 
\\ \\
For proving directly the corresponding ``divergent" Ramanujan-like series for $1/\pi^2$, namely
\[
\sum_{n=0}^{\infty} \frac{\left(\frac12\right)_n^3\left(\frac13\right)_n\left(\frac23\right)_n}{(1)_n^5} (28n^2+18n+3) \,  (-27)^n \, \, `` \! = \! " \, \, \frac{6}{\pi^2}.
\]
we can use \, \texttt{entry:=1/2+2*n,0+3/2*n,1/2+3/2*n,1/2+n,1/2+n;}.

\end{document}